\documentclass[12pt,reqno]{amsart}
\usepackage{amsmath, amssymb, amsfonts, amsthm, xcolor}
\usepackage{mathtools}
\usepackage{hyperref}
\newtheorem{thm}{Theorem}[section]

\newtheorem{prop}[thm]{Proposition}
\theoremstyle{definition}
\newtheorem{dfn}[thm]{Definition}
\newtheorem{exple}[thm]{Example}

\newtheorem{remark}[thm]{Remark}
\theoremstyle{plain}
\newtheorem{cor}[thm]{Corollary}

\numberwithin{equation}{section}
\usepackage[a4paper, top = 1.2in, bottom = 1.2in, left = 1.2in, right = 1.2in]{geometry}
\numberwithin{equation}{section}

\newcommand{\N}{\mathbb{N}}
\newcommand{\Q}{\mathbb{Q}}

\newcommand{\Z}{\mathbb{Z}}

\newcommand{\F}{\mathbb{F}}

\newcommand{\m}{\mathfrak{m}}
\newcommand{\n}{\mathfrak{n}}
\newcommand{\p}{\mathfrak{p}}

\newcommand{\SL}{\mathrm{SL}}
\newcommand{\GL}{\mathrm{GL}}

\newcommand{\ord}{\mathrm{ord}}

\newcommand{\mrm}[1]{\mathrm{#1}}

\def\1{1\!\!1}

\newcommand{\pmat}[4]{ \begin{pmatrix} #1 & #2 \\ #3 & #4 \end{pmatrix}}

\newcommand{\psmat}[4]{\bigl( \begin{smallmatrix} #1 & #2 \\ #3 & #4 \end{smallmatrix} \bigr)}

\def\dis{\displaystyle}

\title[Weakly holomorphic Drinfeld modular forms of level $T$]{A Basis for the space of weakly holomorphic Drinfeld modular forms of level $T$}
\author[T. Dalal]{Tarun Dalal}
\email{tarun.dalal80@gmail.com}
\address{
Department of Mathematics \\
Harish-Chandra Research Institute, HBNI, Chhatnag Road, Jhunsi, Prayagraj 211019, India.
}

\keywords{Drinfeld modular forms, Weakly holomorphic, Generating function, Theta operator}
\subjclass[2010]{11F52}
\date{\today}
\begin{document}
\allowdisplaybreaks
\begin{abstract}

In this article, we explicitly construct a canonical basis for the space of certain weakly holomorphic Drinfeld modular forms for $\Gamma_0(T)$ (resp., for $\Gamma_0^+(T)$) and compute the generating function satisfied by the basis elements. We also give an explicit expression for the action of the $\Theta$-operator, which depends on the divisor of meromorphic Drinfeld modular forms.   
\end{abstract}

\maketitle

\section{Introduction and Statements of the main results}
In the theory of classical modular forms, the existence of a canonical basis for the space of weakly holomorphic modular forms is well known. For example, in \cite{DJ08}, Duke and Jenkins explicitly constructed a canonical basis for the space of weakly holomorphic modular forms for $\SL_2(\Z)$ using the weight $4$ Eisenstein series and the modular discriminant function. Moreover, they computed the generating function satisfied by the basis elements. The generating function plays an important role in studying the zeros of the basis elements. These results have been generalized for other genus zero groups $\Gamma_0(p)$ and $\Gamma_0^+(p)$ by many authors (cf. \cite{Gui09}, \cite{GJ13}, \cite{HJ14} and  \cite{CK13} for more details). 

On the other hand, the existence of a canonical basis for the space of weakly holomorphic Drinfeld modular forms is known only for the case $\GL_2(A)$. Recently in \cite{Cho22}, Choi explicitly constructed a canonical basis for the space of weakly holomorphic Drinfeld modular forms for $\GL_2(A)$. She also computed the generating function satisfied by the basis elements.

In the first part of this article, we explicitly construct a canonical basis for the space of certain weakly holomorphic Drinfeld modular forms for $\Gamma_0(T)$ and compute the generating function satisfied by the basis elements. In the second part of this article, we study similar problems for the space of weakly holomorphic Drinfeld modular forms for $\Gamma^+_0(T)$.

\subsubsection{Notations}
Throughout the article, we fix to use the following notations.

Let $p$ be an odd prime and $q=p^r$ for some $r \in \N$.
Let $\F_q$ denote the finite field of order $q$. We set $A :=\F_q[T]$ 
and $K :=\F_q(T)$.
Let $K_\infty=\F_q((\frac{1}{T}))$ be the completion of $K$ 
with respect to the infinite place $\infty$ (corresponding to $\frac{1}{T}$-adic valuation) 
and the completion of an algebraic closure of $K_{\infty}$ is denoted by $C$. We define the congruence subgroup $\Gamma_0(T)$ of $\GL_2(A)$ 
by $\{\psmat{a}{b}{c}{d}\in \mathrm{GL}_2(A) : c\in (T) \}$.

Let $L=\tilde{\pi}A \subseteq C$ be the $A$-lattice of rank $1$,
corresponding to the rank $1$ Drinfeld module given by $\rho_T=TX+X^q$, which is also known as the Carlitz module, where $\tilde{\pi}\in K_\infty(\sqrt[q-1]{-T})$ is defined up to a $(q-1)$-th root of unity. The Drinfeld upper half-plane $\Omega= C-K_\infty$ has a rigid analytic structure and
the group $\GL_2(K_\infty)$ acts on $\Omega$ via fractional linear transformations.


We now recall the definition of weakly holomorphic Drinfeld modular form for $\Gamma_0(T)$ and discuss some of its basic properties.

\begin{dfn}
\label{Definition of modular forms for level T}
Let $k\in \Z$ and $l\in \Z/(q-1)\Z$.
A meromorphic (resp., holomorphic) function $f: \Omega \rightarrow C$ is said to be a meromorphic (resp., weakly holomorphic) Drinfeld modular form of weight $k$, type $l$ for $\Gamma_0(T)$ if
\begin{enumerate}
\item $f(\gamma z)=(\det\gamma)^{-l} (cz+d)^k f(z)$ for all $\gamma = \pmat{a}{b}{c}{d}\in \Gamma_0(T)$,
\item $f$ is meromorphic at the cusps $\infty$ and $0$.
\end{enumerate}
\end{dfn}
The space of weakly holomorphic Drinfeld modular forms of weight $k$, type $l$ for $\Gamma_0(T)$ will be denoted by $M_{k,l}^!(\Gamma_0(T))$.
Let $M_{k,l}^\#(\Gamma_0(T))$ denote the space of weakly holomorphic Drinfeld modular forms of weight $k$, type $l$ for $\Gamma_0(T)$ which are holomorphic on $\Omega \cup \{0\}$. A weakly holomorphic Drinfeld modular form $f$ of weight $k$, type $l$ for $\Gamma_0(T)$ is said to be a Drinfeld modular form of weight $k$, type $l$ for $\Gamma_0(T)$ if $k\geq 0$ and $f$ is holomorphic at the cusps $\infty$ and $0$. The set of all Drinfeld modular forms of weight $k$, type $l$ for $\Gamma_0(T)$ will be denoted by $M_{k,l}(\Gamma_0(T))$. Observe that, if $k\not \equiv 2l \pmod {q-1}$ and $f$ satisfies condition (1) of Definition \ref{Definition of modular forms for level T}, then $f=0$. 

 Henceforth, we always assume that $k\in \Z$ and $l\in \Z/(q-1)\Z$ such that $k\equiv 2l \pmod {q-1}$.
      Let $0 \leq l \leq q-2$ be a lift of $l \in \Z/(q-1)\Z$. By abuse of notation,
      we continue to write $l$ for the integer as well as its class. Define $r_{k, l} := \frac{k-2l}{q-1}$. Note that, since $q$ is odd, the necessary condition $k\equiv 2l \pmod {q-1}$ implies $k$ is even.

It is well known that the modular curve $X_0(T)=\overline{\Gamma_0(T)\char`\\ \Omega}$ has two cusps which are denoted by $\infty$ and $0$. Moreover the operators $W_T:=\psmat{0}{-1}{T}{0}$ and $\psmat{0}{-1}{1}{0}$ permute the cusps. In order to discuss the expansions of Drinfeld modular forms at the cusps, we introduce the notion of ``slash operator".

Any $x\in K_\infty^\times$ has the unique expression 
$x= \zeta_x\big(\frac{1}{T} \big)^{v_\infty(x)}u_x,$
where $\zeta_x\in \F_q^\times$, and $v_\infty(u_x-1)\geq 0$ ($v_\infty$ is the valuation at $\infty$). For $\gamma = \psmat{a}{b}{c}{d}\in \GL_2(K_{\infty})$ and $f:\Omega\rightarrow C$,
we define 
$$f|_{k,l} \gamma := \zeta_{\det\gamma}^l\big(\frac{\det \gamma}{\zeta_{\det(\gamma)}} \big)^{k/2}(cz+d)^{-k}f(\gamma z).$$ 
By definition, for any $\gamma\in \Gamma_0(T)$, we have $f|_{k,l} \gamma =(\det\gamma)^{l} (cz+d)^{-k} f(\gamma z)$.
\subsubsection{Expansions at the cusps}
Let $u(z) :=  \frac{1}{e_L(\tilde{\pi}z)}$, where $e_L(z):= z{\prod_{\substack{0 \ne \lambda \in L }}}(1-\frac{z}{\lambda}) $
is the exponential function attached to the lattice $L$.
 It is well known that if $f$ is a meromorphic Drinfeld modular form of weight $k$, type $l$ for $\Gamma_0(T)$, then it has $u$-expansion at $\infty$ of the form $\sum_{ i=n}^\infty a_f((q-1)i+l)u^{(q-1)i+l}$, for some $n\in \Z$ and $a_f((q-1)n+l)\ne 0$. We write $\ord_\infty f:=(q-1)n+l$, which is also called the order of vanishing of $f$ at the cusp $\infty$.

On the other hand, if $f$ is a meromorphic Drinfeld modular form of weight $k$, type $l$ for $\Gamma_0(T)$, then it is easy to check that $f|_{k,l}W_T$ is also a meromorphic Drinfeld modular form of weight $k$, type $l$ for $\Gamma_0(T)$ and
\begin{equation}
f|_{k,l}\psmat{0}{-1}{1}{0}(z)=z^{-k}f(\psmat{0}{-1}{1}{0}\frac{z}{T})=T^{-\frac{k}{2}}(f|_{k,l}W_T)(\frac{z}{T}).
\end{equation}
Hence $f|_{k,l}\psmat{0}{-1}{1}{0}$ has a series expansion in $u_0:=u(z/T)$
of the form $\sum_{ i=m}^\infty b_f((q-1)i+l)u_0^{(q-1)i+l}$, for some $m\in \Z$ and $b_f((q-1)m+l)\ne 0$. We write $\ord_0 f:=(q-1)m+l$, which is also called the order of vanishing of $f$ at the cusp $0$.

 We now introduce some important Drinfeld modular forms.
\begin{exple}[\cite{Gos80}, \cite{Gek88}]
\label{Eisenstein Series}
Let $d\in \N$. For $z\in \Omega$, the function
\begin{equation*}
g_d(z) := (-1)^{d+1}\tilde{\pi}^{1-q^d}L_d \sum_{\substack{a,b\in A \\ (a,b)\ne (0,0)}} \frac{1}{(az+b)^{q^d-1}}
\end{equation*}
is a Drinfeld modular form of weight $q^d-1$, type $0$ for $\mathrm{GL}_2(A)$,
where $\tilde{\pi}$ is the Carlitz period and $L_d:=(T^q-T)\cdots(T^{q^d}-T)$ is the least common multiple of all monic polynomials of degree $d$.
\end{exple}
\begin{exple}\cite[\S 4]{DK2}
  Recall that the functions
\begin{equation}
\label{delta_T and delta_W}
 \Delta_T(z) := \frac{g_1(Tz)-g_1(z)}{T^q-T}\ \mrm{and} \ 
\Delta_W(z) := \frac{T^qg_1(Tz)-Tg_1(z)}{T^q-T}
 \end{equation}
both are Drinfeld modular forms of weight $q-1$, type $0$ for $\Gamma_0(T)$ and their $u$-expansions at $\infty$ are given by
$$\Delta_T=u^{q-1}-u^{q(q-1)}+\cdots\in A[[u]], \  \mrm{and}$$
$$\Delta_W = 1+ Tu^{q-1} - T^qu^{q(q-1)} + \cdots \in A[[u]].$$
By definition, we have $\Delta_W|W_T=-T^{\frac{q-1}{2}}\Delta_T$ and $\Delta_T|W_T= -T^{-\frac{q-1}{2}}\Delta_W$. Moreover, it is well known that the modular form $\Delta_T$ (resp., $\Delta_W$) vanishes $q-1$ times at  $\infty$ (resp., at  $0$) and non-zero on $\Omega \cup \{0\}$ (resp.,  on $\Omega \cup \{\infty\}$).
\end{exple} 
\begin{exple}
In \cite{Gek88}, Gekeler defined the function
\begin{equation*}
 E(z):= \frac{1}{\tilde{\pi}} \sum_{\substack{a\in \F_q[T] \\ a \ \mathrm{monic}}} \bigg( \sum_{b\in \F_q[T]} \frac{a}{az+b} \bigg)
\end{equation*}
which is analogous to the Eisenstein series of weight $2$ over $\Q$. 
The function $E$ is not modular, but it satisfies the following transformation rule
\begin{equation}
\label{Etransformation}
E(\gamma z) = (\mathrm{det} \gamma)^{-1} (cz+d)^2E(z) - c\tilde{\pi}^{-1}(\mathrm{det}\gamma)^{-1}(cz+d)
\end{equation}
for $\gamma = \psmat{a}{b}{c}{d}\in \GL_2(A)$. Using the function $E$, we can construct the function $E_T(z):= E(z)-TE(Tz)$ which is a Drinfeld modular form of weight $2$, type $1$ for $\Gamma_0(T)$ (cf. \cite[Proposition 3.3]{DK1}). Since the $u$-expansion of $E$ is $u+ u^{(q-1)^2+1}+\cdots$ and $u(Tz)=u^q+\cdots$, the $u$-expansion of $E_T(z)$ is given by $u-Tu^q+\cdots \in A[[u]]$. Furthermore, we have the relation $E^{q-1}_T=\Delta_W\Delta_T$ (cf. \cite[Proposition 4.3]{DK2}). As a result, $E_T$ vanishes exactly once at the cusps $\infty, 0$ and non-zero elsewhere. Moreover, we have $E_T|W_T=-E_T$ (cf. \cite[Proposition 3.3]{DK1}).
\end{exple}
\subsection{Statements of the main results for $\Gamma_0(T)$}
Our first result is the following:
\begin{thm}
\label{basis for level T intro}
For every integer $i\geq 0$, there exists a unique function $f_{r_{k,l},i}\in M_{k,l}^\#(\Gamma_0(T))$ with $u$-expansion of the form
$$f_{r_{k,l},i}(z)= u^{(q-1)(r_{k,l}-i)+l}+\mathcal{O}(u^{(q-1)(r_{k,l}+1)+l}).$$
Furthermore, the coefficients of the $u$-expansion of $f_{r_{k,l},i}$ all belong to $A$ and there exists a unique monic polynomial $F_{r_{k,l},i}(x)\in A[x]$ of degree $i$ such that $f_{r_{k,l},i}=F_{r_{k,l},i}(j_T) \Delta_T^{r_{k,l}}E_T^l$, where $j_T:=\frac{\Delta_W}{\Delta_T}$.

Consequently, the set $\{f_{r_{k,l},i}: i \geq 0\}$ forms a basis for $M_{k,l}^\#(\Gamma_0(T))$. 
\end{thm}
Let $C_{r_{k,l},i}$ be the constant term of the polynomial $F_{r_{k,l},i}(x)$. 
We prove that the basis elements $f_{r_{k,l},i}$'s and the constant terms $C_{r_{k,l},i}$'s satisfy the following relations:
\begin{thm}
\label{functional equation level T intro}
We have
\begin{enumerate}
\item[(i)] $\sum_{i\geq 0} f_{r_{k,l},i}(\tau)u^{(q-1)(-r_{k,l}+i)+1-l}=\frac{E_T^q f_{r_{k,l},0}(\tau)}{\big( j_T(z)-j_T(\tau)\big)\Delta_T^2 f_{r_{k,l},0}(z)}, \ \mathrm{for} \ \tau \in \Omega.$
\item [(ii)] $\sum_{i\geq 0}C_{r_{k,l},i}u^{(q-1)(-r_{k,l}+i)+1-l}=\Delta_T^{-r_{k,l}}E_T^{1-l}.$
\end{enumerate}
 \end{thm}
Observe that for $k=l=0$, from Theorem \ref{functional equation level T intro} (ii) we have $E_T=\sum_{i\geq 0} C_{r_{0,0},i}u^{(q-1)i+1}$. This gives a new characterization of the function $E_T$.

The proof of Theorem \ref{functional equation level T intro} heavily depends on the action of the $\Theta$-operator on meromorphic Drinfeld modular forms. The $\Theta$-operator is defined as 
$$\Theta := \frac{1}{\tilde{\pi}}\frac{d}{dz}= -u^2\frac{d}{du}.$$

We prove the following explicit expression for the action of the $\Theta$-operator.
\begin{thm}
\label{Theta f /f combine form Intro}
If $f$ is a meromorphic Drinfeld modular form of weight $k$, type $l$ for $\Gamma_0(T)$, then
$$\Theta(f)=-kEf+ f \cdot\Big(\sum_{\substack{\pi(\tau)\in \Gamma_0(T)\char`\\ \Omega}} \ord_\tau f \cdot \frac{E_T^q(z)}{\big(j_T(\tau)- j_T(z)\big)\Delta_T^2(z)} + \ord_0f\cdot E_T\Big),
$$
where $\pi : \Omega \rightarrow \Gamma_0(T)\char`\\ \Omega$ is the quotient map and  $\ord_\tau f$ denotes the mod-$p$ reduction of the order of vanishing of $f$ at $\tau$.
\end{thm}

In the second part of the article, we study Drinfeld modular forms for $\Gamma_0^+(T) :=\langle \Gamma_0(T), W_T\rangle$. 
In \cite{Cho09a}, Choi defined meromorphic Drinfeld modular forms for $\Gamma_0^+(T)$. Due to the choice of normalization, the author restricted herself to the case $k=2l$.
We now propose an alternative definition of meromorphic Drinfeld modular forms for $\Gamma_0^+(T)$ so that the assumption $k=2l$ is redundant.

\begin{dfn}
A meromorphic (resp., holomorphic) function $f: \Omega \rightarrow C$ is said to be a meromorphic (resp., weakly holomorphic) Drinfeld modular form of weight $k$, type $l$ for $\Gamma_0^+(T)$ if
\begin{enumerate}
\item $f|_{k,l} \gamma = f$ for all $\gamma\in \Gamma_0^+(T)$,
\item $f$ is meromorphic at the cusp $\infty$.
\end{enumerate}
\end{dfn}
The space of weakly holomorphic Drinfeld modular forms of weight $k$, type $l$ for $\Gamma_0^+(T)$ will be denoted by $M_{k,l}^!(\Gamma_0^+(T))$. 
A weakly holomorphic Drinfeld modular form $f$ of weight $k$, type $l$ for $\Gamma_0^+(T)$ is said to be a Drinfeld modular form of weight $k$, type $l$ for $\Gamma_0^+(T)$ if $k\geq 0$ and $f$ is holomorphic at the cusp $\infty$. The set of all Drinfeld modular forms of weight $k$, type $l$ for $\Gamma_0^+(T)$ will be denoted by $M_{k,l}(\Gamma_0^+(T))$.

Every meromorphic Drinfeld modular form $f$ of weight $k$, type $l$ for $\Gamma_0^+(T)$ has $u$-expansion at $\infty$ of the form $\sum_{ i=n}^\infty a_f((q-1)i+l)u^{(q-1)i+l}$, for some $n\in \Z$ and $a_f((q-1)n+l)\ne 0$. We define $\ord_\infty(f):=(q-1)n+l$.

For $k$ and $l$ as before we define
$d^+:=
\begin{cases}
\frac{r_{k,l}+1}{2},  &\mathrm{if\ } r_{k,l} \ \mathrm{is\ odd}\\
\frac{r_{k,l}}{2},  &\mathrm{if\ } r_{k,l} \ \mathrm{is\ even\ and\ } l \ \mathrm{is\ odd}\\
\frac{r_{k,l}}{2}+1,  &\mathrm{if\ } r_{k,l} \ \mathrm{is\ even\ and\ } l \ \mathrm{is\ even}.
\end{cases}
$

\subsection{Statements of the main results for $\Gamma_0^+(T)$}
For simplicity of notations, we write $j_T^+:=\frac{(\Delta_W-T^{\frac{q-1}{2}}\Delta_T)^2}{E_T^{q-1}}$.
We first construct the following canonical basis for $M_{k,l}^!(\Gamma_0^+(T))$. 
\begin{thm}
\label{basis for + space intro}
For every integer $i\geq 0$, there exists a unique function $f_{r_{k,l},i}^+\in M_{k,l}^!(\Gamma_0^+(T))$ with $u$-expansion of the form
$$f_{r_{k,l},i}^+(z)= u^{(q-1)((d^+-1)-i)+l}+\mathcal{O}(u^{(q-1)d^++l}).$$
Furthermore, the coefficients of the $u$-expansion of $f^+_{r_{k,l},i}$ all belong to $A$ and there exists a unique monic polynomial $F_{r_{k,l},i}^+(x)\in A[x]$ of degree $i$ such that 
\begin{enumerate}
\item $f_{r_{k,l},i}^+=F_{r_{k,l},i}^+(j_T^+) \big(\Delta_W - (-1)^{l}T^{\frac{q-1}{2}} \Delta_T\big)E_T^{(d^+-1)(q-1)+l}$, when $r_{k,l}$ is odd.
\item {$f_{r_{k,l},i}^+=F_{r_{k,l},i}^+(j_T^+) \big(\Delta_W^2 -T^{q-1} \Delta_T^2\big)E_T^{(d^+-1)(q-1)+l}$}, when $r_{k,l}$ is even and $l$ is odd.
\item $f_{r_{k,l},i}^+=F_{r_{k,l},i}^+(j_T^+)E_T^{(d^+-1)(q-1)+l}$, when $r_{k,l}$ is even and $l$ is even.
\end{enumerate}

Consequently, the set $\{f^+_{r_{k,l},i}: i \geq 0\}$ forms a basis for $M_{k,l}^!(\Gamma_0^+(T))$.
\end{thm}
We also prove that the basis elements $f^+_{r_{k,l},i}$'s satisfy the following relation:

\begin{thm}
\label{functional equation for + space}
For any $\tau \in \Omega$, we have
$$\sum_{i\geq 0} f^+_{r_{k,l},i}(\tau)u^{(q-1)(-(d^+-1)+i)+1-l} =\frac{(\Delta_W^2-T^{q-1}\Delta_T^2)f^+_{r_{k,l},0}(\tau)}{E_T^{q-2}\big(j^+_T(\tau)- j^+_T(z)\big)f^+_{r_{k,l},0}(z)}.$$
\end{thm}
{As in the previous case, the proof of Theorem \ref{functional equation for + space} heavily depends on the action of the $\Theta$-operator on meromorphic Drinfeld modular forms for $\Gamma_0^+(T)$. We prove the following result regarding the action of the $\Theta$-operator.}
\begin{thm}
\label{value of theta f/f for + intro}
If $f$ is a meromorphic Drinfeld modular form of weight $k$, type $l$ for $\Gamma^+_0(T)$, then
$$\Theta(f)=- \frac{k}{2}(E(z)+TE(Tz))f+ f\cdot\Big(\sum_{\substack{\pi^+(\tau)\in \Gamma^+_0(T)\char`\\ \Omega}}\frac{\ord_\tau f}{e^+_\tau} \cdot \frac{\Delta_W^2-T^{q-1}\Delta_T^2}{E_T^{q-2}\big(j_T^+(\tau)-j_T^+(z)\big)}\Big),
$$
where $\pi^+ : \Omega \rightarrow \Gamma^+_0(T)\char`\\ \Omega$ is the quotient map and $\ord_\tau f$ denotes the mod-$p$ reduction of the order of vanishing of $f$ at $\tau$
(cf. \S\ref{A basis for +} for the definition of $e^+_\tau$).
\end{thm}

Since the methods of this article work for any degree one prime ideal, it is interesting to study similar results for $\Gamma_0(\n)$ with $\deg \n>1$.

\section{ A Basis for $M_{k,l}^\#(\Gamma_0(T))$}

In this section, we prove Theorem \ref{basis for level T intro}, Theorem \ref{functional equation level T intro} and Theorem \ref{Theta f /f combine form Intro}.
We start by giving an upper bound on the order of vanishing at $\infty$.

\begin{prop}
\label{order of vanishing at infty}
For $f\ne 0 \in M_{k,l}^\#(\Gamma_0(T))$, we have that $\ord_\infty f \leq k-l$.
\end{prop}
\begin{proof}
If $\ord_\infty f > k-l$, then $G:=\frac{f}{E_T^l}\Delta_T^{-r_{k,l}}$ is holomorphic on $\Omega\cup \{0,\infty\}$ and $G\in M_{0,0}(\Gamma_0(T))$. Hence we have 
$$G=\mathrm{constant} \ne 0 , \ i.e. \ f\Delta_T^{-r_{k,l}} =\mathrm{constant} \cdot E_T^l.$$
 Which contradicts that $E_T$ vanishes exactly once at the cusp $\infty$. This proves the result.
\end{proof}
Note that $j_T:= \frac{\Delta_W}{\Delta_T}\in M_{0,0}^\#(\Gamma_0(T))$ has the $u$-expansion $u^{-(q-1)}+T+\mathcal{O}(u^{q-1})$.
We are now ready to prove Theorem \ref{basis for level T intro}. For the benefit of the readers, we recall the statement of Theorem \ref{basis for level T intro}. 
\begin{thm}
\label{existence of Victor Miller basis for level T}
For every integer $i\geq 0$, there exists a unique function $f_{r_{k,l},i}\in M_{k,l}^\#(\Gamma_0(T))$ with $u$-expansion of the form
$$f_{r_{k,l},i}(z)= u^{(q-1)(r_{k,l}-i)+l}+\mathcal{O}(u^{(q-1)(r_{k,l}+1)+l}).$$
Furthermore, the coefficients of the $u$-expansion of $f_{r_{k,l},i}$ all belong to $A$ and there exists a unique monic polynomial $F_{r_{k,l},i}(x)\in A[x]$ of degree $i$ such that $f_{r_{k,l},i}=F_{r_{k,l},i}(j_T) \Delta_T^{r_{k,l}}E_T^l$.

Consequently, the set $\{f_{r_{k,l},i}: i \geq 0\}$ forms a basis for $M_{k,l}^\#(\Gamma_0(T))$. 
\end{thm}
\begin{proof}
We choose $f_{r_{k,l},0}=\Delta_T^{r_{k,l}}E_T^l$. For any $i\in \N$, we can construct $f_{r_{k,l},i}$ from $j_T^i \Delta_T^{r_{k,l}}E_T^l$ by subtracting off suitable constant multiples of $f_{r_{k,l},i-d}$ for $1\leq d \leq i$ (observe that all those constants belong to $A$). The uniqueness of $f_{r_{k,l},i}$'s follows from Proposition \ref{order of vanishing at infty}. Consequently, there exists a unique monic polynomial $F_{r_{k,l},i}(x)\in A[x]$ of degree $i$ such that $f_{r_{k,l},i}=F_{r_{k,l},i}(j_T) \Delta_T^{r_{k,l}}E_T^l$ and the set $\{f_{r_{k,l},i}: i \geq 0\}$ forms a basis for $M_{k,l}^\#(\Gamma_0(T))$. Since the coefficients of the $u$-expansions of $f_{r_{k,l},0}$ and $j_T$ all belong to $A$, by construction the coefficients of the $u$-expansion of $f_{r_{k,l},i}$ all belong to $A$.
\end{proof}
\begin{remark}
Observe that, if $k=l=0$, then by construction we have $f_{r_{0,0},0}=1, f_{r_{0,0},1}=\frac{\Delta_W}{\Delta_T}-T$. 
For simplicity of notations, we write $f_{0,i}:=f_{r_{0,0},i}$ and $C_{0,i}:=C_{r_{0,0},i}$ for $i\geq 0$.
\end{remark}

{For $\tau\in \Omega$, we define $e_\tau:= [\Gamma_0(T)_\tau/(\Gamma_0(T)_\tau \cap Z(K))]$, where $\Gamma_0(T)_\tau$ is the stabilizer of $\tau$ in $\Gamma_0(T)$ and $Z(K)$ is the set of scalar matrices. It is well known that $e_\tau$ is either $1$ or $q+1$. Thus $e_\tau \equiv 1 \pmod p$.}
\begin{prop}
\label{formulas for residues}
Let $G$ be a meromorphic Drinfeld modular form of weight $2$, type $1$ for $\Gamma_0(T)$ and $\pi : \Omega \rightarrow \Gamma_0(T)\char`\\ \Omega$ be the quotient map.
\begin{enumerate}
\item [(i)] \label{residue P1}
For $\tau \in \Omega$, we have $\mathrm{Res}_{\pi(\tau)}G(z)dz=\frac{1}{e_\tau}\mathrm{Res}_\tau G(z)=\mathrm{Res}_\tau G(z)$. 
\item [(ii)] \label{residue P2}
If the $u$-expansion of $G$ at $\infty$ is given by
      $G(z)=\sum_{i\geq -n} a_{G,\infty} (i(q-1)+1)u^{i(q-1)+1},$ then 
      $\mathrm{Res}_\infty G(z)dz=\frac{a_{G,\infty} (1)}{\tilde{\pi}}$,
      where $\tilde{\pi}$ is the Carlitz period.
\item [(iii)] \label{residue P3}
If the $u_0$-expansion of $G|_{2,1}\psmat{0}{-1}{1}{0}$ is given by $\sum_{i\geq -n} a_{G,0}(i(q-1)+1)u_0^{i(q-1)+1}$ where $u_0:=u(z/T)$,
      then $\mathrm{Res}_0\ G(z)dz=\frac{T}{\tilde{\pi}}a_{G,0}(1)$.
\end{enumerate}
\end{prop}
\begin{proof}
The proof of (i) is analogous to \cite[Lemma 3.1(ii)]{Cho06}. The proofs of (ii) and (iii) can be found in \cite[Proposition 3.3]{DK3}.
\end{proof}

It is well known that if $f$ is a meromorphic Drinfeld modular form of weight $k$, type $l$ for $\Gamma_0(T)$, then $\partial_k f:= \Theta(f)+kEf$ is a meromorphic Drinfeld modular form of weight $k+2$, type $l+1$ for $\Gamma_0(T)$. Consequently, $\frac{\Theta(f)}{f}+kE$ is a meromorpic Drinfeld modular form of weight $2$, type $1$ for $\Gamma_0(T)$. 

\begin{thm}
\label{theta f / f general value for level T}
Let $f$ be a non-zero meromorphic Drinfeld modular form of weight $k$, type $l$ for $\Gamma_0(T)$ with $u$-expansion
$f(z)= u^{(q-1)n+l}+ \sum_{i=n+1}^\infty a_f((q-1)i+l)u^{(q-1)i+l}.$
Then 
$$\frac{\Theta(f)}{f}=-kE+  \sum_{\pi(\tau)\in \Gamma_0(T)\char`\\ \Omega} \mathrm{ord}_\tau f\Big(-u- \sum_{i=1}^\infty f_{0,i}(\tau) u^{(q-1)i+1} \Big) +\ord_0 f\Big(u+ \sum_{i=0}^\infty C_{0,i}u^{(q-1)i+l} \Big),$$
where $C_{0,i}$ is the constant term of the polynomial $F_{r_{0,0},i}(x)$.

\end{thm}
\begin{proof}
For any $i\in \N \cup \{0\}$, let $G_i(z)=(\frac{\Theta(f)}{f}+kE)f_{0,i}(z)$, where $f_{0,0}(z)=1$. Then $G_i(z)$ is a meromorphic Drinfeld modular form of weight $2$, type $1$ for $\Gamma_0(T)$. Note that we can write 
\begin{equation}
\label{u expansion of theta f divided by f}
\frac{\Theta(f)}{f}= -\sum_{i=0}^\infty b_{(q-1)i}u^{(q-1)i+1}
\end{equation}
where $b_0:= (q-1)n+l$ and $b_{(q-1)i}\in C$ for $i\in \N$. Now
\begin{align*}
G_i(z)&= \Big(\frac{\Theta(f)}{f}+k E\Big)f_{0,i}(z)\\
&=\Big(-\sum_{i=0}^\infty b_{(q-1)i}u^{(q-1)i+1}+k\sum_{a\in A_+}au(az)\Big)\cdot (u^{-i(q-1)}+\mathcal{O}(u^{(q-1)}))\\
&= \cdots + \Big(-b_{(q-1)i} + k\sum_{a\in A_+}a\{u_a\}_{(q-1)i+1}\Big)u+\cdots,
\end{align*}
where $\{u_a\}_{(q-1)i+1}$ denotes the coefficient of $u^{(q-1)i+1}$ in $u(az)$.

By Proposition \ref{formulas for residues}, we get
\begin{equation}
\mathrm{Res}_\infty G_i(z)dz= \frac{1}{\tilde{\pi}}\Big(-b_{(q-1)i} +k \sum_{a\in A_+}a\cdot \{u_a\}_{(q-1)i+1} \Big).
\end{equation}
Since $E(z)$ and $f_{0,i}(z)$ both are holomorphic on $\Omega$, a simple computation shows that
\begin{equation}
\mathrm{Res}_{\pi(\tau)} G_i(z)dz=\mathrm{Res}_\tau G(z)= \frac{\mathrm{ord}_\tau f \cdot f_{0,i}(\tau)}{\tilde{\pi}}.
\end{equation}
We now compute $\mrm{Res}_0 G_i(z)dz$. Let 
$\gamma:=\psmat{0}{-1}{1}{0}$ and $ g:= f|_{k,l} \gamma$, i.e.  $g(z)=z^{-k}f(\gamma z)$. 
Since $\frac{dz}{d(\gamma z)}=z^2$, we have
\begin{align*}
\frac{\Theta(f)}{f}(\gamma z)+ kE(\gamma z)= \frac{z^2}{\tilde{\pi}f(\gamma z)} \frac{d}{dz}(f(\gamma z)) + kE(\gamma z)&= \frac{z^2}{\tilde{\pi}z^kg(z)} \frac{d}{dz}(z^kg(z)) + kE(\gamma z) \\
&=z^2 \frac{\Theta(g)}{g} + k z^2E(z) \ (\mathrm{cf}. \eqref{Etransformation}).
\end{align*}
Thus
\begin{equation}
\Big(\frac{\Theta(f)}{f}+ kE \Big)|_{2,1}\pmat{0}{-1}{1}{0}= \frac{\Theta (g)}{g} +kE.
\end{equation}

Let the $u_0$-expansion of $g$ be of the form $\sum_{j\geq m}a_{g,0}(j(q-1)+l)u_0^{j(q-1)+l}$ ($m\in \Z$) and the $u_0$-expansion of $\frac{1}{g}$ be of the form $\sum_{j\geq -m}a^\prime(j(q-1)-l)u_0^{(j(q-1))-l}$. Since $dz= -\frac{T}{\tilde{\pi}}u_0^{-2}du_0$, the $u_0$-expansion of $\frac{\Theta(g)}{g}$ is given by
$$-\frac{1}{T}(j(q-1)+l)u_0+\mathcal{O}(u_0^q)=-\frac{\ord_0 f}{T}u_0+\mathcal{O}(u_0^q).$$
Since the $u$-expansion of $E$ is of the form $u+\cdots$ and $u(z)=\frac{u_0^q}{1+Tu_0^{q-1}}=u_0^q+\cdots$, the $u_0$-expansion of $E$ is of the form $u_0^q+\cdots$. 

Recall that $C_{0,i}$ is the constant term of the polynomial $F_{r_{0,0},i}(x)$ for $i\in \N\cup \{0\}$ and $C_{0,0}=1$.
By Theorem \ref{existence of Victor Miller basis for level T}, we can write $f_{0,i}=(\frac{\Delta_W}{\Delta_T})^i+c_{0,1}(\frac{\Delta_W}{\Delta_T})^{i-1}+\cdots+c_{0,i-1}\frac{\Delta_W}{\Delta_T} + C_{0,i}$. 
{Since $\Delta_W|W_T=-T^{\frac{q-1}{2}}\Delta_T$ and $\Delta_T|W_T= -T^{-\frac{q-1}{2}}\Delta_W$,}
the $u_0$-expansion of $f_{0,i}|_{0,0}\psmat{0}{-1}{1}{0}$ is of the form $C_{0,i}+ \mathcal{O}(u_0^{q-1})$.

Hence we get the $u_0$-expansion
\begin{align*}
\Big(\frac{\Theta(f)}{f}+ kE \Big)f_{0,i}|_{2,1}\psmat{0}{-1}{1}{0}&= \Big(-\frac{\ord_0 f}{T}u_0+\mathcal{O}(u_0^q)\Big)\cdot \Big( C_{0,i}+ \mathcal{O}(u_0^{q-1})\Big)= -\frac{\ord_0 f}{T}C_{0,i}u_0+\mathcal{O}(u_0^q).
\end{align*}
By Proposition \ref{formulas for residues}, we have $\mathrm{Res}_0 G_i(z)dz=-\frac{\ord_0 f}{\tilde{\pi}}C_{0,i}$. 

Finally, by residue theorem we get 
\begin{equation}
\label{coefficients relation using residue theorem}
b_{(q-1)i}= k \sum_{a\in A_+}a\cdot \{u_a\}_{(q-1)i+1}  + \sum_{\pi(\tau)\in \Gamma_0(T)\char`\\ \Omega} \mathrm{ord}_\tau f \cdot f_{0,i}(\tau) -{\ord_0 f} \cdot C_{0,i}.
\end{equation}
Since $u_a=u^{q^{\deg(a)}}+\cdots$, from \eqref{coefficients relation using residue theorem} we obtain
\begin{equation}
\label{value of b0 for level T}
b_0=(q-1)n+l=-\frac{(q-1)n+l}{q-1}= k+ \sum_{\pi(\tau)\in \Gamma_0(T)\char`\\ \Omega} \mathrm{ord}_\tau f  -{\ord_0 f},
\end{equation}
which can be rewritten as
\begin{equation}
\label{Valence formula for level T equation}
\sum_{\pi(\tau)\in \Gamma_0(T)\char`\\ \Omega}\frac{\ord_\tau(f)}{e_\tau}+ \frac{\ord_\infty(f)}{q-1} +\frac{\ord_0 f}{q-1}= \frac{k}{q-1}
\end{equation}
(Note that the we are considering the above sum over $\F_p$).

Since $E=\sum_{i= 0}^\infty \Big(\sum_{a\in A_+}a\{u_a\}_{(q-1)i+1}\Big)u^{(q-1)i+1}$, using \eqref{coefficients relation using residue theorem} we get
{

\begin{align*}
&\frac{\Theta(f)}{f}=-\sum_{i=0}^\infty b_{(q-1)i}u^{(q-1)i+1}\\
&= -k\sum_{i=0}^\infty \Big( \sum_{a\in A_+}a\cdot \{u_a\}_{(q-1)i+1}\Big) u^{(q-1)i+1} - \sum_{i=0}^\infty\Big( \sum_{\pi(\tau)\in \Gamma_0(T)\char`\\ \Omega} \mathrm{ord}_\tau f \cdot f_{0,i}(\tau)\Big)u^{(q-1)i+1}\\
& \hspace{12cm} +{\ord_0 f}\sum_{i=0}^\infty C_{0,i}u^{(q-1)i+l}.
\end{align*}}
Hence we conclude that 
\begin{equation}
\label{value of Theta f / f for level T}
\frac{\Theta(f)}{f}=-kE+  \sum_{\pi(\tau)\in \Gamma_0(T)\char`\\ \Omega} \mathrm{ord}_\tau f\Big(-u- \sum_{i=1}^\infty f_{0,i}(\tau) u^{(q-1)i+1} \Big) +\ord_0 f\Big(u+ \sum_{i=1}^\infty C_{0,i}u^{(q-1)i+l} \Big).
\end{equation}
\end{proof}

We define $\Phi_T(z,0):=\frac{\Theta(j_T(z))}{j_T(z)}$ and $\Phi_T(z,\tau):=\frac{\Theta(j_T(z)-j_T(\tau))}{j_T(z)-j_T(\tau)}$ for fixed $\tau \in \Omega$. Then the functions $\Phi_T(z,0)$ and $\Phi(z, \tau)$ both are meromorphic Drinfeld modular forms of weight $2$, type $1$ for $\Gamma_0(T)$. 

Note that if $\ord_x(j_T(z)-j_T(\tau))>0$ for $x\in \Omega$, then $j_T(x)=j_T(\tau)$, and consequently $f_{0,i}(x)=f_{0,i}(\tau)$. Therefore by
 \eqref{value of Theta f / f for level T} and \eqref{value of b0 for level T} the $u$-expansions of $\Phi_T(z,0)$ and $\Phi_T(z,\tau)$ are given by 
$$\Phi_T(z,0) = -u- \sum_{i= 1}^\infty C_{0,i}u^{(q-1)i+1} \ \mathrm{and} \ \Phi_T(z,\tau) = -u- \sum_{i= 1}^\infty f_{0,i}(\tau)u^{(q-1)i+1}.$$

Thus we obtain the following result.

\begin{prop}
 The function $\Phi_T(z,0) := -u- \sum_{i= 1}^\infty C_{0,i}u^{(q-1)i+1}$ and for each $\tau\in \Omega$ the function
$\Phi_T(z,\tau):= -u- \sum_{i= 1}^\infty f_{0,i}(\tau)u^{(q-1)i+1}$
are all meromorphic Drinfeld modular forms of weight $2$, type $1$ for $\Gamma_0(T)$.
\end{prop}
We now compute the explicit forms of the functions $\Phi_T(z,0)$ and $\Phi_T(z,\tau)$.
\begin{prop}
\label{value of phi T}
We have 
$\Phi_T(z,0)= -\frac{E_T^q(z)}{j_T(z)\Delta_T^2(z)}=-E_T$, and
$\Phi_T(z,\tau)=\frac{E_T^q(z)}{\big(j_T(\tau)- j_T(z)\big)\Delta_T^2(z)}$ for $\tau\in \Omega$.
\end{prop}

\begin{proof}
Since $\partial_{q-1}\Delta_T(z)=0, \partial_{q-1}(\Delta_W)=-\Delta_WE_T$ (cf. \cite[Proposition 4.3]{DK2}) and $j_T(z)=\frac{\Delta_W(z)}{\Delta_T(z)}$, we can write 
\begin{align*}
\Delta_W\frac{\Theta(j_T(z))}{j_T(z)} &= \Delta_W\frac{\Theta(j_T(z))}{j_T(z)} + j_T(z) \partial_{q-1}\Delta_T(z)\\
&= \Delta_T \Theta(j_T(z)) + j_T(z)\partial_{q-1}\Delta_T(z)\\
&= \partial_{q-1}(\Delta_T(z)j_T(z))=-\Delta_WE_T.
\end{align*}
Therefore $\Theta(j_T(z))=-E_T\frac{\Delta_W(z)}{\Delta_T(z)}$. Now, the result follows from the equality $E_T^{q-1}=\Delta_W\Delta_T$ (cf.~\cite[Proposition 4.3]{DK2}).
\end{proof}

 Comparing the $u$-expansions of $\dis \Phi_T(z,0)$ and $E_T=\sum_{\substack{a\in A_+\\ T\nmid a}}au(az)=u+\sum_{\substack{a\in A_+\\ T\nmid a, , a\neq 1}}au(az)$, by Proposition \ref{value of phi T} we get the following corollary:
\begin{cor}
For $i\in \N$ we have
$C_{0,i}=\sum_{\substack{a\in A_+\\ T\nmid a, , a\neq 1}}a\{u_a\}_{(q-1)i+1}$.
\end{cor}

Finally, combining Proposition \ref{value of phi T} and Theorem \ref{theta f / f general value for level T} we obtain the following theorem:
\begin{thm}
\label{Theta f /f combine form}
If $f$ is a meromorphic Drinfeld modular form of weight $k$, type $l$ for $\Gamma_0(T)$, then
$$\Theta(f)=
-kEf+ f\cdot\Big(\sum_{\substack{\pi(\tau)\in \Gamma_0(T)\char`\\ \Omega}} \ord_\tau f \cdot \frac{E_T^q(z)}{\big(j_T(\tau)- j_T(z)\big)\Delta_T^2(z)} + \ord_0f\cdot E_T\Big).
$$
\end{thm}

We now compute the generating function satisfied by the functions $f_{r_{k,l}}$'s.

Let the $u$-expansions of $f_{r_{k,l},0}$ and $\frac{1}{f_{r_{k,l},0}}$ be given by
$${f_{r_{k,l},0}(z)= \sum_{i\geq 0}a_{f_{r_{k,l},0}}((q-1)(r_{k,l}+i)+l)u^{(q-1)(r_{k,l}+i)+l} \ \mathrm{and}} $$
$$\frac{1}{f_{r_{k,l},0}(z)}= \sum_{i\geq 0}a^\prime_{f_{r_{k,l},0}}((q-1)(-r_{k,l}+i)-l)u^{(q-1)(-r_{k,l}+i)-l},$$
where $a_{f_{r_{k,l},0}}((q-1)r_{k,l}+l)=a^\prime_{f_{r_{k,l},0}}(-(q-1)r_{k,l}-l)=1$.

\begin{prop}
\label{required relation between fi and f0}
For any $i\in \N$ we have
\begin{enumerate}
\item [(i)] $f_{r_{k,l},i}(z)= f_{r_{k,l},0}(z)\sum_{n+s=i}a^\prime_{f_{r_{k,l},0}}((q-1)(-r_{k,l}+n)-l)f_{0,s}(z),$
\item [(ii)]$C_{r_{k,l},i}= \sum_{n+s=i}a^\prime_{f_{r_{k,l},0}}((q-1)(-r_{k,l}+n)-l)C_{0,s}$.
\end{enumerate}
\end{prop}
\begin{proof}
By the relation $f_{r_{k,l},0}\cdot \frac{1}{f_{r_{k,l},0}}=1$ and the corresponding $u$-expansions we get
\begin{equation}
\label{equation obtain after considering f.finverse =1}
\sum_{j=0}^r a_{f_{r_{k,l},0}}((q-1)(r_{k,l}+j)+l)a^\prime_{f_{r_{k,l},0}}((q-1)(-r_{k,l}+r-j)-l)=0, \ \mathrm{for} \ r\geq 1.
\end{equation}
On the other hand
\begin{align*}
&f_{r_{k,l},0}(z)\sum_{n+s=i}a^\prime_{f_{r_{k,l},0}}((q-1)(-r_{k,l}+n)-l)f_{0,s}(z)\\
&= \sum_{j\geq 0}a_{f_{r_{k,l},0}}((q-1)(r_{k,l}+j)+l)u^{(q-1)(r_{k,l}+j)+l} \sum_{n=0}^i a^\prime_{f_{r_{k,l},0}}((q-1)(-r_{k,l}+n)-l) f_{0,i-n}(z)\\
&= \sum_{j\geq 0}a_{f_{r_{k,l},0}}((q-1)(r_{k,l}+j)+l)u^{(q-1)(r_{k,l}+j)+l}\\
& \hspace{5cm} \cdot \sum_{n=0}^i a^\prime_{f_{r_{k,l},0}}((q-1)(-r_{k,l}+n)-l) (u^{-(i-n)(q-1)}+\mathcal{O}(u^{q-1}))\\
&= u^{(q-1)(r_{k,l}-i)+l}\\
&\hspace{2cm} + \sum_{r=1}^i \Big(\sum_{j+n=r} a_{f_{r_{k,l},0}}((q-1)(r_{k,l}+j)+l)a^\prime_{f_{r_{k,l},0}}((q-1)(-r_{k,l}+n)-l)\Big)u^{(q-1)(r_{k,l}-i+r)+l}\\
&\hspace{10cm} + \mathcal{O}(u^{(q-1)(r_{k,l}+1)+l})\\
&= u^{(q-1)(r_{k,l}-i)+l}+ \mathcal{O}(u^{(q-1)(r_{k,l}+1)+l}) \ (\mathrm{by} \ \eqref{equation obtain after considering f.finverse =1}).
\end{align*}
By the uniqueness of $f_{r_{k,l},i}$'s (cf. Theorem \ref{existence of Victor Miller basis for level T}), we conclude that 
\begin{equation}
\label{required relation for level T equation}
f_{r_{k,l},i}(z)= f_{r_{k,l},0}(z)\sum_{n+s=i}a^\prime_{f_{r_{k,l},0}}((q-1)(-r_{k,l}+n)-l)f_{0,s}(z).
\end{equation}
This proves (i).

{
Recall that $C_{r_{k,l},i}$ is the constant term of the polynomial $F_{r_{k,l},i}(x)$ and $C_{r_{k,l},0}=1.$ It is easy to check that the $u_0$-expansion of $f_{r_{k,l},i}|_{k,l}\psmat{0}{-1}{1}{0}$ is of the form $C_{r_{k,l},i}(\frac{-1}{T^{q-1}})^{r_{k,l}}(\frac{-1}{T})^lu_0^l +\mathcal{O}(u_0^{(q-1)+l})$.} Considering the $u_0$-expansions, from \eqref{required relation for level T equation} we have
\begin{align*}
C_{r_{k,l},i}(-\frac{1}{T^{q-1}})^{r_{k,l}}(-\frac{1}{T})^lu^l +\mathcal{O}(u_0^{(q-1)+l})&=  ((-\frac{1}{T^{q-1}})^{r_{k,l}}(-\frac{1}{T})^lu_0^l +\mathcal{O}(u_0^{(q-1)+l}))\\ & \hspace{1cm} \cdot \sum_{n+s=i}a^\prime_{f_{r_{k,l},0}}((q-1)(-r_{k,l}+n)-l)(C_{0,s}+\mathcal{O}(u_0^{q-1})).
\end{align*}
Comparing the coefficients of $u_0^l$ we get the equality
\begin{equation}
C_{r_{k,l},i}= \sum_{n+s=i}a^\prime_{f_{r_{k,l},0}}((q-1)(-r_{k,l}+n)-l)C_{0,s}.
\end{equation}
This proves (ii).
\end{proof}

Finally, we are ready to proof Theorem \ref{functional equation level T intro}.  
\begin{thm}
\label{functional relation among basis elements}
The basis elements $f_{r_{k,l},i}$'s and the constants $C_{r_{k,l},i}$'s satisfy the following relations:
\begin{enumerate}
\item[(i)] $\sum_{i\geq 0} f_{r_{k,l},i}(\tau)u^{(q-1)(-r_{k,l}+i)+1-l}=\frac{E_T^q f_{r_{k,l},0}(\tau)}{\big( j_T(z)-j_T(\tau)\big)\Delta_T^2 f_{r_{k,l},0}(z)}, \ \mathrm{for} \ \tau \in \Omega.$
\item [(ii)] $\sum_{i\geq 0}C_{r_{k,l},i}u^{(q-1)(-r_{k,l}+i)+1-l}=\Delta_T^{-r_{k,l}}E_T^{1-l}.$
\end{enumerate}
\end{thm}
\begin{proof}
By Proposition \ref{value of phi T}, we have
\begin{equation}
\label{required equation between proving the relation between fi's}
\frac{E_T^q f_{r_{k,l},0}(\tau)}{\big( j_T(z)-j_T(\tau)\big)\Delta_T^2 f_{r_{k,l},0}(z)}= -\frac{f_{r_{k,l},0}(\tau)}{f_{r_{k,l},0}(z)} \Phi_T(z,\tau).
\end{equation}
On the other hand, considering the $u$-expansions we obtain
\begin{align*}
&\frac{f_{r_{k,l},0}(\tau)}{f_{r_{k,l},0}(z)} \Phi_T(z,\tau)\\
&= -f_{r_{k,l},0}(\tau) \sum_{i\geq 0}a^\prime_{f_{r_{k,l},0}}((q-1)(-r_{k,l}+i)-l)u^{(q-1)(-r_{k,l}+i)-l} \sum_{i\geq 0}f_{0,i}(\tau)u^{(q-1)i+1} \\
&= - \sum_{i\geq 0}f_{r_{k,l},0}(\tau)\Big(\sum_{r+s=i} a^\prime_{f_{r_{k,l},0}}((q-1)(-r_{k,l}+r)-l)f_{0,s}(\tau) \Big) u^{(q-1)(-r_{k,l}+i)+1-l}\\
&= -\sum_{i\geq 0} f_{r_{k,l},i}(\tau) u^{(q-1)(-r_{k,l}+i)+1-l} \ (\mathrm{cf. \ Proposition} \ \ref{required relation between fi and f0}). 
\end{align*}
Now (i) follows from \eqref{required equation between proving the relation between fi's}.

Using Proposition \ref{required relation between fi and f0}, a similar calculation as above gives
$$-\sum_{i\geq 0}C_{r_{k,l},i}u^{(q-1)(-r_{k,l}+i)+1-l}=\frac{\Phi_T(z,0)}{f_{r_{k,l},0}(z)}=-\Delta_T^{-r_{k,l}}E_T^{1-l}.$$
This proves (ii).
\end{proof}

\begin{remark}
Recall that $\Delta_W|W_T=-T^{\frac{q-1}{2}}\Delta_T, \Delta_T|W_T= -T^{-\frac{q-1}{2}}\Delta_W$ and the operator $W_T$ permutes the cusps $0$ and $\infty$. By interchanging the roles of $\Delta_W$ and $\Delta_T$, and using the similar arguments discussed so far, we can construct a canonical basis for the space of weakly holomorphic Drinfeld modular forms which are holomorphic on $\Omega\cup\{\infty\}$. Furthermore, we can compute the generating function satisfied by the basis elements. 
\end{remark}

\section{A basis for $M_{k,l}(\Gamma_0^+(T))$}
 Recall that $W_T:=\psmat{0}{-1}{T}{0}$, and $\Gamma_0^+(T):=\langle \Gamma_0(T), W_T\rangle$.
 Since $W_T$ permutes the cusps $\infty$ and $0$, the modular curve $X_0^+(T)=\overline{\Gamma_0^+(T)\char`\\ \Omega}$ has only cusp, we call it the cusp at $\infty$. In order to prove Theorem \ref{basis for + space intro}, we first need to construct a basis for $M_{k,l}(\Gamma_0^+(T))$  ($k\geq 0$).

 Recall that $\Delta_T|W_T=-T^{-\frac{q-1}{2}}\Delta_W$, $\Delta_W|W_T=-T^{\frac{q-1}{2}}\Delta_T$ and the set $\{\Delta_W^{r_{k,l}-i}\Delta_T^iE_T^l: 0\leq i \leq r_{k,l}\}$ forms a basis for $M_{k,l}(\Gamma_0(T))$ ($k\geq 0$) (cf. \cite[Proposition 4.3]{DK2}). 
 
For $k\geq 0$, we define the set $M_{k,l}(\Gamma_0(T))^+:=\{f\in M_{k,l}(\Gamma_0(T)): f|_{k,l}W_T=f\}$.
 
 Since the operator $W_T$ permutes the cusps $\infty$ and $0$ of $\Gamma_0(T)$, for $k\geq 0$ we have the equality $M_{k,l}(\Gamma_0^+(T))=M_{k,l}(\Gamma_0(T))^+$.

Observe that for $0\leq i \leq \lfloor \frac{r_{k,l}-1}{2} \rfloor$, we have
\begin{equation}
\Big( \Delta_W^{r_{k,l}-i}\Delta_T^i + (-1)^{l+r_{k,l}}T^{\frac{q-1}{2}r_{k,l}-i(q-1)}\Delta_W^i\Delta_T^{r_{k,l}-i} \Big)E_T^l\in M_{k,l}(\Gamma_0(T))^+.
\end{equation}
Consider the set $$S:=\Big\{ \Big( \Delta_W^{r_{k,l}-i}\Delta_T^i + (-1)^{l+r_{k,l}}T^{\frac{q-1}{2}r_{k,l}-i(q-1)}\Delta_W^i\Delta_T^{r_{k,l}-i} \Big)E_T^l: 0\leq i \leq \lfloor \frac{r_{k,l}-1}{2}\rfloor \Big\}.$$
\begin{thm}
\label{basis for holomorphic Drinfeld modular forms for +}
\begin{enumerate}
\item If $r_{k,l}$ is odd, then the set $S$ forms a basis for $M_{k,l}(\Gamma_0(T))^+$.
\item If $r_{k,l}$ is even and $l$ is odd, then the set $S$ forms a basis for $M_{k,l}(\Gamma_0(T))^+$.
\item If $r_{k,l}$ is even and $l$ is even, then the set $S\cup \{\Delta_W^{\frac{r_{k,l}}{2}}\Delta_T^{\frac{r_{k,l}}{2}}E_T^l\}$ forms a basis for $M_{k,l}(\Gamma_0(T))^+$.
\end{enumerate}
Consequently, we obtain 
$$\dim M_{k,l}(\Gamma_0(T))^+=
\begin{cases}
\frac{r_{k,l}+1}{2},  &\mathrm{if\ } r_{k,l} \ \mathrm{is\ odd}\\
\frac{r_{k,l}}{2},  &\mathrm{if\ } r_{k,l} \ \mathrm{is\ even\ and\ } l \ \mathrm{is\ odd}\\
\frac{r_{k,l}}{2}+1,  &\mathrm{if\ } r_{k,l} \ \mathrm{is\ even\ and\ } l \ \mathrm{is\ even}.
\end{cases}
$$
\end{thm}
\begin{proof} 
Clearly, the elements of the set $S\cup \{\Delta_W^{\frac{r_{k,l}}{2}}\Delta_T^{\frac{r_{k,l}}{2}}E_T^l\}$ are linearly independent.
Hence it is enough to prove that the set $S\cup \{\Delta_W^{\frac{r_{k,l}}{2}}\Delta_T^{\frac{r_{k,l}}{2}}E_T^l\}$ generates $M_{k,l}(\Gamma_0(T))^+$.
 Let $f (\ne 0)\in M_{k,l}(\Gamma_0(T))^+$. Since the set $\{\Delta_W^{r_{k,l}-i}\Delta_T^iE_T^l: 0\leq i \leq r_{k,l}\}$ forms a basis for $M_{k,l}(\Gamma_0(T))$, there exist $c_i\in C$ (not all zero) such that
\begin{equation}
\label{linear combination of f}
f= \sum_{i=0}^{r_{k,l}}c_i \Delta_W^{r_{k,l}-i}\Delta_T^iE_T^l.
\end{equation}
Therefore 
\begin{equation}
\label{linear combination of f|WT}
f|_{k,l}W_T=\sum_{i=0}^{r_{k,l}}c_i \Delta_W^{r_{k,l}-i}\Delta_T^iE_T^l|_{k,l}W_T= \sum_{i=0}^{r_{k,l}}c_i (-1)^{r_{k,l}+l}T^{\frac{q-1}{2}r_{k,l}-(q-1)i}\Delta_W^i\Delta_T^{r_{k,l}-i}E_T^l.
\end{equation}
Since $f|_{k,l}W_T=f$, from \eqref{linear combination of f} and \eqref{linear combination of f|WT}, we have
\begin{equation*}
\sum_{i=0}^{r_{k,l}}c_i \Delta_W^{r_{k,l}-i}\Delta_T^iE_T^l = \sum_{i=0}^{r_{k,l}}c_i (-1)^{r_{k,l}+l}T^{\frac{q-1}{2}r_{k,l}-(q-1)i}\Delta_W^i\Delta_T^{r_{k,l}-i}E_T^l.
\end{equation*}
Adding the terms with equal powers we get
\begin{equation}
\label{relation after adding terms}
\sum_{i=0}^{r_{k,l}}(c_i-(-1)^{r_{k,l}+l}c_{r_{k,l}-i}T^{-\frac{q-1}{2}r_{k,l}+(q-1)i})\Delta_W^{r_{k,l}-i}\Delta_T^i=0.
\end{equation}
Since $\Delta_W$ and $\Delta_T$ are algebraically independent, from \eqref{relation after adding terms} we obtain
\begin{equation}
\label{relation between ci's for basis of + space}
c_{r_{k,l}-i}=(-1)^{r_{k,l}+l}T^{\frac{q-1}{2}r_{k,l}-(q-1)i}c_{i}, \ \mathrm{for} \ 0\leq i \leq r_{k,l}.
\end{equation}
When $r_{k,l}$ is odd, using \eqref{relation between ci's for basis of + space} we can write
\begin{equation*}
f= \sum_{i=0}^{\frac{r_{k,l}-1}{2}}c_i \Big(\Delta_W^{r_{k,l}-i}\Delta_T^i + (-1)^{l+r_{k,l}}T^{\frac{q-1}{2}r_{k,l}-i(q-1)}\Delta_W^i\Delta_T^{r_{k,l}-i} \Big)E_T^l.
\end{equation*}
When $r_{k,l}$ is even and $l$ is odd, using \eqref{relation between ci's for basis of + space} we can write
\begin{equation*}
f= \sum_{i=0}^{\frac{r_{k,l}}{2}-1}c_i \Big(\Delta_W^{r_{k,l}-i}\Delta_T^i + (-1)^{l+r_{k,l}}T^{\frac{q-1}{2}r_{k,l}-i(q-1)}\Delta_W^i\Delta_T^{r_{k,l}-i} \Big)E_T^l.
\end{equation*}
When $r_{k,l}$ is even and $l$ is even, using \eqref{relation between ci's for basis of + space} we can write
\begin{equation*}
f= \sum_{i=0}^{\frac{r_{k,l}}{2}-1}c_i \Big(\Delta_W^{r_{k,l}-i}\Delta_T^i + (-1)^{l+r_{k,l}}T^{\frac{q-1}{2}r_{k,l}-i(q-1)}\Delta_W^i\Delta_T^{r_{k,l}-i} \Big)E_T^l + c_{\frac{r_{k,l}}{2}}\Delta_W^{\frac{r_{k,l}}{2}}\Delta_T^{\frac{r_{k,l}}{2}}E_T^l.
\end{equation*}
The result follows.
\end{proof}

As an immediate consequence of Theorem \ref{basis for holomorphic Drinfeld modular forms for +} we obtain:
 \begin{cor}
\label{order of vanishing at infty of holomorphic modular form for + space}
For every non-zero $f\in M_{k,l}(\Gamma_0(T))^+= M_{k,l}(\Gamma_0^+(T))$, we have $$\ord_\infty f\leq (\dim M_{k,l}(\Gamma_0(T))^+-1)(q-1)+l.$$
\end{cor}
\begin{proof}
We give the complete proof only for the case when $r_{k,l}$ is odd. The proofs for the other cases are analogous to this.
Let $f (\ne 0) \in M_{k,l}(\Gamma_0^+(T))$ and $r_{k,l}$ is odd. Then there exists constants $c_i\in C$ such that 
\begin{equation*}
f= \sum_{i=0}^{\frac{r_{k,l}-1}{2}} c_i\Big( \Delta_W^{r_{k,l}-i}\Delta_T^i + (-1)^{l+r_{k,l}}T^{\frac{q-1}{2}r_{k,l}-i(q-1)}\Delta_W^i\Delta_T^{r_{k,l}-i} \Big)E_T^l.
\end{equation*}
The result follows since the $u$-expansions of the basis elements are given by

$\Big( \Delta_W^{r_{k,l}-i}\Delta_T^i + (-1)^{l+r_{k,l}}T^{\frac{q-1}{2}r_{k,l}-i(q-1)}\Delta_W^i\Delta_T^{r_{k,l}-i} \Big)E_T^l=u^{i(q-1)+l}+\cdots.$
\end{proof}

\section{A basis for $M_{k,l}^!(\Gamma_0^+(T))$}
\label{A basis for +}
In this section, we prove Theorem \ref{basis for + space intro}, Theorem \ref{functional equation for + space} and Theorem \ref{value of theta f/f for + intro}. 
We start by giving an upper bound on the order of vanishing at $\infty$.

For integers $k,l$ with $0\leq l < q-1$ and $k\equiv 2l \pmod {q-1}$ we define
$$d^+:=
\begin{cases}
\frac{r_{k,l}+1}{2},  &\mathrm{if\ } r_{k,l} \ \mathrm{is\ odd}\\
\frac{r_{k,l}}{2},  &\mathrm{if\ } r_{k,l} \ \mathrm{is\ even\ and\ } l \ \mathrm{is\ odd}\\
\frac{r_{k,l}}{2}+1,  &\mathrm{if\ } r_{k,l} \ \mathrm{is\ even\ and\ } l \ \mathrm{is\ even}.
\end{cases}
$$
\begin{prop}
\label{order of vanishing at infty of weakly holomorphic DMF of levle T + space}
For every non zero $f\in M_{k,l}^!(\Gamma_0^+(T))$, we have $$\ord_\infty f\leq (d^+-1)(q-1)+l.$$
\end{prop}
\begin{proof}
Let $f\ne 0 \in M_{k,l}^!(\Gamma_0^+(T))$ such that $\ord_\infty f> (d^+-1)(q-1)+l$. Then $f E_T^{-(d^+-1)(q-1)}\in M_{(r_{k,l}-2(d^+-1))(q-1)+2l,l}(\Gamma_0^+(T))$ and 
\begin{equation}
\label{contradiction for vanishing at infty of holomorphic DMF + space}
\ord_\infty(f E_T^{-(d^+-1)(q-1)})> l.
\end{equation}
 On the other hand, using Corollary \ref{order of vanishing at infty of holomorphic modular form for + space}, a simple calculation shows that $\ord_\infty(f E_T^{-(d^+-1)(q-1)})\leq l.$ Which contradicts \eqref{contradiction for vanishing at infty of holomorphic DMF + space}. This proves the proposition.
\end{proof}
Note that $j_T^+:=\frac{(\Delta_W-T^{\frac{q-1}{2}}\Delta_T)^2}{E_T^{q-1}}=\frac{1}{u^{q-1}}+(2(T-T^{\frac{q-1}{2}})-T)+\mathcal{O}(u^{q-1})\in M_{0,0}^!(\Gamma_0^+(T))$.

Now we are ready to prove Theorem \ref{basis for + space intro}.
\begin{thm}
For every integer $i\geq 0$, there exists a unique function $f_{r_{k,l},i}^+\in M_{k,l}^!(\Gamma_0^+(T))$ with $u$-expansion of the form
$$f_{r_{k,l},i}^+(z)= u^{(q-1)((d^+-1)-i)+l}+\mathcal{O}(u^{(q-1)d^++l}).$$
Furthermore, the coefficients of the $u$-expansions of $f_{r_{k,l},i}$ all belong to $A$ and there exists a unique monic polynomial $F_{r_{k,l},i}^+(x)\in A[x]$ of degree $i$ such that 
\begin{enumerate}
\item $f_{r_{k,l},i}^+=F_{r_{k,l},i}^+(j_T^+) \big(\Delta_W - (-1)^{l}T^{\frac{q-1}{2}} \Delta_T\big)E_T^{(d^+-1)(q-1)+l}$, when $r_{k,l}$ is odd.
\item {$f_{r_{k,l},i}^+=F_{r_{k,l},i}^+(j_T^+) \big(\Delta_W^2 -T^{q-1} \Delta_T^2\big)E_T^{(d^+-1)(q-1)+l}$}, when $r_{k,l}$ is even and $l$ is odd.
\item $f_{r_{k,l},i}^+=F_{r_{k,l},i}^+(j_T^+)E_T^{(d^+-1)(q-1)+l}$, when $r_{k,l}$ is even and $l$ is even.
\end{enumerate}

Consequently, the set $\{f_{r_{k,l},i}: i \geq 0\}$ forms a basis for $M_{k,l}^!(\Gamma_0^+(T))$.
\end{thm}
\begin{proof}
We give the complete proof only for the case when $r_{k,l}$ is odd. The proofs for the other cases are analogous to this. Let $r_{k,l}$ be odd.

We choose $f^+_{r_{k,l},0}=\big(\Delta_W - (-1)^{l}T^{\frac{q-1}{2}} \Delta_T\big)E_T^{(d^+-1)(q-1)+l}$. Furthermore, for any $i\in \N$ we can construct $f^+_{r_{k,l},i}$ from $\big(j_T^+\big)^i \big(\Delta_W - (-1)^{l}T^{\frac{q-1}{2}} \Delta_T\big)E_T^{(d^+-1)(q-1)+l}$ by subtracting off suitable constant multiples of $f^+_{r_{k,l},i-j}$ for $1\leq j \leq i$ (observe that all those constants belong to $A$). The uniqueness follows from Proposition \ref{order of vanishing at infty of weakly holomorphic DMF of levle T + space}. Consequently, there exists a unique monic polynomial $F^+_{r_{k,l},i}(x)\in A[x]$ of degree $i$ such that 
$f_{r_{k,l},i}^+=F_{r_{k,l},i}^+(j_T^+) \big(\Delta_W - (-1)^{l}T^{\frac{q-1}{2}} \Delta_T\big)E_T^{(d^+-1)(q-1)+l},$
and the set $\{f^+_{r_{k,l},i}: i \geq 0\}$ forms a basis for $M_{k,l}^!(\Gamma_0^+(T))$. Since the coefficients of the $u$-expansions of $f^+_{r_{k,l},0}$ and $j^+_T$ all belong to $A$, by construction the coefficients of the $u$-expansion of $f^+_{r_{k,l},i}$ all belong to $A$.
\end{proof}
For simplicity of notation, we write $f^+_{0,i}:=f^+_{r_{0,0},i}$ for $i\geq 0$. By construction, we have $f^+_{0,0}=1$.


It is easy to check that if $G(z)$ is a meromorphic Drinfeld modular form of weight $2$, type $1$ for $\Gamma_0^+(T)$, then $G(z)dz$ is a $1$-form on $X_0^+(T)$. Since $W_T$ does not fix the cusp $\infty$, $u^{q-1}$ acts as the parameter at $\infty$ for $X_0^+(T)$.

{For $\tau\in \Omega$, we define $e^+_\tau:= [\Gamma^+_0(T)_\tau/(\Gamma^+_0(T)_\tau \cap Z(K))]$, where $\Gamma^+_0(T)_\tau$ is the stabilizer of $\tau$ in $\Gamma^+_0(T)$ and $Z(K)$ is the set of scalar matrices. It is well known that $e^+_\tau$ is either $e_\tau$ or $2e_\tau$. Thus $e^+_\tau \equiv 1,2 \pmod p$.} We recall the following result:
\begin{prop}\cite[Lemma 3.2]{Cho09a}
\label{residue formula for +}
Let $G$ be a meromorphic Drinfeld modular form of weight $2$, type $1$ for $\Gamma^+_0(T)$ and $\pi^+ : \Omega \rightarrow \Gamma^+_0(T)\char`\\ \Omega$ be the quotient map. 
\begin{enumerate}
\item For $\tau \in \Omega$ we have $\mathrm{Res}_{\pi^+(\tau)}G(z)dz=\frac{1}{e^+_\tau}\mathrm{Res}_\tau G(z)$.
\item If the $u$-expansion of $G$ at $\infty$ is given by
      $G(z)=\sum_{i\geq -n} a_{G,\infty} (i(q-1)+1)u^{i(q-1)+1},$ then 
      $\mathrm{Res}_\infty G(z)dz=\frac{a_{G,\infty} (1)}{\tilde{\pi}}$.
\end{enumerate}
\end{prop}

An easy computation shows that if $f$ is a meromorphic Drinfeld modular form of weight $k$, type $l$ for $\Gamma_0^+(T)$, then
\begin{equation}
\partial^+(f):= \big(\Theta f +\frac{k}{2}(E(z)+TE(Tz))f\big)= \partial (f) -\frac{k}{2}E_Tf
\end{equation}
 is a meromorphic Drinfeld modular form of weight $k+2$, type $l+1$ for $\Gamma_0^+(T)$ (cf. \cite[Proposition 3.4]{DK1} for more detail)
and 
$\big(\frac{\Theta f}{f} +\frac{k}{2}(E(z)+TE(Tz))\big)$ is a meromorphic Drinfeld modular form of weight $2$, type $1$ for $\Gamma_0^+(T)$. Moreover for $i=1,2$, if $f_i$ is a meromorphic Drinfeld modular form of weight $k_i$, type $l_i$ for $\Gamma_0(T^+)$, then $\partial^+(f_1f_2)=f_1\partial^+(f_2)+f_2\partial^+(f_1).$

\begin{thm}
\label{theta f/f general value for +}
Let $f$ be a meromorphic Drinfeld modular form of weight $2$, type $1$ for $\Gamma_0^+(T)$ with $u$-expansion
$f(z)= u^{(q-1)n+l}+ \sum_{i=n+1}^\infty a_f((q-1)i+l)u^{(q-1)i+l}.$ Then
$$\frac{\Theta(f)}{f} + \frac{k}{2}(E(z)+TE(Tz))= \sum_{\pi^+(\tau) \in \Gamma_0^+(T)\char`\\\Omega} \frac{\ord_\tau f}{e^+_\tau}\Big(-u- \sum_{i\geq 1} f^+_{0,i}(\tau)u^{(q-1)i+1}\Big).$$
\end{thm}
\begin{proof}
For any $i\in \N\cup \{0\}$, let $G_i(z)=(\frac{\Theta(f)}{f}+\frac{k}{2} (E(z)+TE(Tz)))f^+_{0,i}(z)$ where $f^+_{0,0}=1$. Then $G_i(z)$ is a meromorphic Drinfeld modular form of weight $2$, type $1$ for $\Gamma_0^+(T)$. As before we can write 
\begin{equation}
\label{u expansion of theta f divided by f for +}
\frac{\Theta(f)}{f}= -\sum_{i=0}^\infty b_{(q-1)i}u^{(q-1)i+1},
\end{equation}
where $b_0=(q-1)n+l=-n+l$ and $b_{(q-1)i}\in C$ for $i\in \N$. Now
\begin{align*}
G_i(z)&= \Big(\frac{\Theta(f)}{f}+\frac{k}{2} (E(z)+TE(Tz))\Big)f^+_{0,i}(z)\\
&=\Big(-\sum_{i=0}^\infty b_{(q-1)i}u^{(q-1)i+1}+\frac{k}{2}(\sum_{a\in A_+}au(az)+T\sum_{a\in A_+}au(Taz))\Big)\\
&\hspace{8cm} \cdot (u^{-i(q-1)}+\mathcal{O}(u^{(q-1)}))\\
&= \cdots + \Big(-b_{(q-1)i} + \frac{k}{2}(\sum_{a\in A_+}a\{u_a\}_{(q-1)i+1}+\sum_{a\in A_+}Ta\{u_{Ta}\}_{(q-1)i+1})\Big)u+\cdots,
\end{align*}
where $\{u_a\}_{(q-1)i+1}$ denotes the coefficient of $u^{(q-1)i+1}$ in $u(az)$.

Thus, from Proposition \ref{residue formula for +}, we obtain $$\mathrm{Res}_\infty G_i(z)dz=\frac{1}{\tilde{\pi}}\Big({-b_{(q-1)i} + \frac{k}{2}(\sum_{a\in A_+}a\{u_a\}_{(q-1)i+l}+\sum_{a\in A_+}Ta\{u_{Ta}\}_{(q-1)i+l})}\Big).$$
Since $E(z), E(Tz)$ and $f^+_{0,i}(z)$ are holomorphic on $\Omega$, for $\tau \in \Omega$ we have
$$\mathrm{Res}_\tau G_i(z)= \mathrm{Res}_\tau \frac{\Theta (f)}{f}f^+_{0,i}(z)=\frac{\ord_\tau f}{\tilde{\pi}} \cdot f^+_{0,i}(\tau).$$

Consequently, for $\tau\in \Omega$, Proposition \ref{residue formula for +} we gives $\mathrm{Res}_{\pi^+(\tau)}G(z)dz=\frac{\ord_\tau f}{\tilde{\pi}e^+_{\tau}} \cdot f^+_{0,i}(\tau)$.

Therefore, by the residue theorem, we get
$$\frac{1}{\tilde{\pi}}\Big({-b_{(q-1)i} + \frac{k}{2}(\sum_{a\in A_+}a\{u_a\}_{(q-1)i+l}+\sum_{a\in A_+}Ta\{u_{Ta}\}_{(q-1)i+l})}\Big)+ \sum_{\pi^+(\tau) \in \Gamma_0^+(T)\char`\\\Omega} \frac{\ord_\tau f}{\tilde{\pi}e^+_\tau} \cdot f^+_{0,i}(\tau)=0.$$
Hence 
\begin{equation}
\label{coefficients relation using residue theorem for + }
{b_{(q-1)i} = \frac{k}{2}\Big(\sum_{a\in A_+}a\{u_a\}_{(q-1)i+l}+\sum_{a\in A_+}Ta\{u_{Ta}\}_{(q-1)i+l}\Big)}+ \sum_{\pi^+(\tau) \in \Gamma_0^+(T)\char`\\\Omega} \frac{\ord_\tau f}{e^+_\tau} \cdot f^+_{0,i}(\tau).
\end{equation}
Since $u_a=u^{q^{\deg(a)}}+\ldots$, from \eqref{coefficients relation using residue theorem for + } we obtain
\begin{equation}
\label{value of b0 for +}
b_0=-n+l=(q-1)n+l=-\frac{(q-1)n+l}{q-1}=\frac{k}{2}+ \sum_{\pi^+(\tau) \in \Gamma_0^+(T)\char`\\\Omega} \frac{\ord_\tau f}{e^+_\tau},
\end{equation}
which can be rewritten as
\begin{equation}
\label{Valence formula for + space equation}
\sum_{\pi^+(\tau) \in \Gamma_0^+(T)\char`\\\Omega} \frac{\ord_\tau f}{e^+_\tau}+ \frac{\ord_\infty(f)}{q-1}=\frac{k}{2(q-1)}.
\end{equation}
(Note that the we are considering the above sum over $\F_p$).

Since $\dis E=u+\sum_{i= 1}^\infty \big(\sum_{a\in A_+}a\{u_a\}_{(q-1)i+1}\big)u^{(q-1)i+1}$ and $\dis TE(Tz)=\sum_{i= 1}^\infty \big(\sum_{a\in A_+}Ta\{u_{Ta}\}_{(q-1)i+1}\big)u^{(q-1)i+1}$, combining \eqref{u expansion of theta f divided by f for +} and \eqref{coefficients relation using residue theorem for + } we get

\begin{align*}
&\frac{\Theta(f)}{f}=-\sum_{i=0}^\infty b_{(q-1)i}u^{(q-1)i+1}=-\frac{k}{2}(E(z)+TE(Tz)) - \sum_{i=0}^\infty\Big( \sum_{\pi^+(\tau) \in \Gamma_0^+(T)\char`\\\Omega} \frac{\ord_\tau f}{e^+_\tau} \cdot f^+_{0,i}(\tau)\Big)u^{(q-1)i+1}.
\end{align*}
Hence we conclude that
\begin{equation}
\label{value of Theta f /f for +}
\frac{\Theta(f)}{f}= \sum_{\pi^+(\tau) \in \Gamma_0^+(T)\char`\\\Omega} \frac{\ord_\tau f}{e^+_\tau}(-u- \sum_{i\geq 1} f^+_{0,i}(\tau)u^{(q-1)i+1}) - \frac{k}{2}(E(z)+TE(Tz)).
\end{equation}
This completes the proof.
\end{proof}

Now fix $\tau \in \Omega$. Observe that if $\ord_x(j^+_T(z)-j^+_T(\tau))>0$ for $x\in \Omega$, then $j^+_T(x)=j^+_T(\tau)$, and consequently $f^+_{0,i}(x)=f^+_{0,i}(\tau)$.

\begin{prop}
\label{weight 2, type 1 function for +}
For each $\tau\in \Omega$, the function
$\Phi^+_T(z,\tau):= -u- \sum_{i\geq 1} f^+_{0,i}(\tau)u^{(q-1)i+1}$
is a meromorphic Drinfeld modular form of weight $2$, type $1$ for $\Gamma_0^+(T)$.
\end{prop}
\begin{proof}
For a fixed $\tau\in \Omega$, let $H^+_\tau(z):=\frac{\Theta(j_T^+(z)-j_T^+(\tau))}{j_T^+(z)-j_T^+(\tau)}$. Then $H_\tau(z)$ is a meromorphic Drinfeld modular form of weight $2$, type $1$ for $\Gamma_0^+(T)$. By \eqref{value of Theta f /f for +} and \eqref{value of b0 for +} the $u$-expansion of $H^+_\tau(z)$ is given by 
$$H^+_\tau(z):= -u- \sum_{i\geq 1} f^+_{0,i}(\tau)u^{(q-1)i+1}.$$
The result follows.
\end{proof}
We now compute the explicit form of the function $\Phi_T^+(z,\tau)$.
\begin{prop}
\label{value of Phi_T for +}
For any $\tau \in \Omega$, we have $\Phi_T^+(z,\tau)=\frac{\Delta_W^2-T^{q-1}\Delta_T^2}{E_T^{q-2}\big(j_T^+(\tau)-j_T^+(z)\big)}$.
\end{prop}
\begin{proof}
From the proof of Proposition \ref{weight 2, type 1 function for +}, we have $\Phi_T^+(z, \tau)=\frac{\Theta(j_T^+(z))}{j_T^+(z)-j_T^+(\tau)}$.

Using \cite[Proposition 4.3]{DK2}, it is easy to check that 
\begin{equation}
\partial^+(E_T^{q-1})=0 \ \mathrm{ and} \ \partial^+(\Delta_W-T^{\frac{q-1}{2}}\Delta_T)= -\frac{1}{2}(\Delta_W+T^{\frac{q-1}{2}}\Delta_T)E_T.
\end{equation}
Recall that $j_T^+=\frac{(\Delta_W-T^{\frac{q-1}{2}}\Delta_T)^2}{E_T^{q-1}}$. Now
\begin{align*}
(\Delta_W-T^{\frac{q-1}{2}}\Delta_T)^2 \frac{\Theta(j^+_T)}{j^+_T}
&= (\Delta_W-T^{\frac{q-1}{2}}\Delta_T)^2 \frac{\Theta(j^+_T)}{j^+_T} + j^+_T\partial^+(E_T^{q-1})\\
&= E_T^{q-1}\Theta(j^+_T) + j^+_T\partial^+(E_T^{q-1})\\
&= \partial^+(j^+_T E_T^{q-1})\\
&= \partial^+((\Delta_W-T^{\frac{q-1}{2}}\Delta_T)^2)\\
&= -(\Delta_W^2-T^{q-1}\Delta_T^2)E_T.
\end{align*}
Therefore $\Theta(j^+_T(z))=-\frac{(\Delta_W^2-T^{q-1}\Delta_T^2)}{E_T^{q-2}}$ and the result follows.
\end{proof}
Combining Proposition \ref{value of Phi_T for +} and Theorem \ref{theta f/f general value for +} we obtain the following theorem:
\begin{thm}
If $f$ is a meromorphic Drinfeld modular form of weight $k$, type $l$ for $\Gamma^+_0(T)$, then
$$\Theta(f)=
 -\frac{k}{2}(E(z)+TE(Tz))f+ f\cdot\Big( \sum_{\substack{\pi^+(\tau)\in \Gamma^+_0(T)\char`\\ \Omega}}\frac{\ord_\tau f}{e^+_\tau} \cdot \frac{\Delta_W^2-T^{q-1}\Delta_T^2}{E_T^{q-2}\big(j_T^+(\tau)-j_T^+(z)\big)}\Big).
$$
\end{thm}

We now compute the generating function satisfied by the elements $f^+_{r_{k,l}}$'s.

Let the $u$-expansions of $f^+_{r_{k,l},0}$ and $\frac{1}{f^+_{r_{k,l},0}}$ be given by
$${f^+_{r_{k,l},0}(z)= \sum_{i\geq 0}a_{f^+_{r_{k,l},0}}((q-1)(d^+-1+i)+l)u^{(q-1)(d^+-1+i)+l} \ \mathrm{and}} $$
$$\frac{1}{f^+_{r_{k,l},0}(z)}= \sum_{i\geq 0}a^\prime_{f^+_{r_{k,l},0}}((q-1)(-(d^+-1)+i)-l)u^{(q-1)(-(d^+-1)+i)-l},$$
where $a_{f^+_{r_{k,l},0}}((q-1)(d^+-1)+l)=a^\prime_{f^+_{r_{k,l},0}}(-(q-1)(d^+-1)-l)=1$.

\begin{prop}
\label{required relation between fi and f0 for +}
For any $i\in \N$ we have
$${f^+_{r_{k,l},i}(z)= f^+_{r_{k,l},0}(z)\sum_{n+s=i}a^\prime_{f^+_{r_{k,l},0}}((q-1)(-(d^+-1)+n)-l)f^+_{0,s}(z)}.$$
\end{prop}
\begin{proof}
The proof is similar to the proof of Proposition \ref{required relation between fi and f0}, hence we omit the detail.
\end{proof}

Finally, we are ready to proof Theorem \ref{functional equation for + space}.
\begin{thm}
\label{generating function for + space in text}
The basis elements $f^+_{r_{k,l},i}$'s satisfy the following relation:
$$\sum_{i\geq 0} f^+_{r_{k,l},i}(\tau)u^{(q-1)(-(d^+-1)+i)+1-l} =\frac{(\Delta_W^2-T^{q-1}\Delta_T^2)f^+_{r_{k,l},0}(\tau)}{E_T^{q-2}\big(j^+_T(z)- j^+_T(\tau)\big)f^+_{r_{k,l},0}(z)}, \ \mathrm{for} \ \tau \in \Omega.$$
\end{thm}
\begin{proof}
The proof is similar to the proof of Theorem \ref{functional relation among basis elements}.
\begin{align*}
&\frac{f^+_{r_{k,l},0}(\tau)}{f^+_{r_{k,l},0}(z)} \Phi^+_T(z,\tau)\\
&=  -{f^+_{r_{k,l},0}(\tau)}\sum_{i\geq 0}a^\prime_{f^+_{r_{k,l},0}}((q-1)(-(d^+-1)+i)-l)u^{(q-1)(-(d^+-1)+i)-l} \sum_{i\geq 0} f^+_{0,i}(\tau)u^{(q-1)i+1}\\
&=- \sum_{i\geq 0}f^+_{r_{k,l},0}(\tau)\Big(\sum_{r+s=i} a^\prime_{f^+_{r_{k,l},0}}((q-1)(-(d^+-1)+r)-l)f^+_{0,s}(\tau) \Big)u^{(q-1)(-(d^+-1)+i)+1-l}\\
&= -\sum_{i\geq 0} f^+_{r_{k,l},i}(\tau)u^{(q-1)(-(d^+-1)+i)+1-l}.
\end{align*}
Now the result follows from Proposition \ref{value of Phi_T for +}.
\end{proof}

\section*{Acknowledgments}
The author thanks Dept. of Atomic Energy, Govt of India for the financial support provided to carry out this research work at Harish-Chandra Research Institute.

\bibliographystyle{plain, abbrv}

\end{document}